\documentclass[12pt]{article}
\usepackage {amssymb}
\usepackage {amsmath}
\begin{document}

\begin{center}
\textbf{Phase Flows and Vector Hamiltonians}
\bigskip

\textbf{V.N.Dumachev}\footnote{E-mail: dumv@comch.ru} \\
Voronezh Institute of the MVD of the Russia
\end{center}

\begin{small}

We present a generalization of the Nambu mechanics on the base of
Liouville's theorem. We prove that the Poisson structure of an
$n$-dimensional multisymplectic phase space is induced by
$(n-1)$-Hamiltonian $k$-vector field seach of which requires
introduction of $k$-Hamiltonians.

\end{small}

\bigskip
\textbf{1.} Consider a system of differential equations as a
submanifold $\Sigma $ in a jet bundle $J^n(\pi)$: $E \to M$, defined
by equations ([1], P.10)
\[
\label{eq1} g(t,x_0,x_1,...,x_n)=0,
\]
where $t \in M \subset R$, $u = x_{0} \in U \subset R$, $x_i \in
J^i(\pi) \subset R^n$, $E=M \times U$.  For equations $\Sigma
\subset J^n(\pi)$ which admit a Poisson structure with bracket
\begin{equation}
\label{eq1}\{H,G\}=X_H\rfloor dG=L_{X_H}G
\end{equation}
the Cartan distribution can be given as follows:
\begin{equation}
\label{eq2} \theta _i=dx_i-\{H,x_i\}dt.
\end{equation}
Here $L_{X_H}$ is the Lie derivative with respect to $X_H \in
\Lambda ^1$, $\Lambda ^n$ is the exterior graded algebra of
$k$-polyvector fields, $H=H(\textbf{x})$ is an unknown function.

\bigskip
\textbf{2.} The classical symplectic mechanics on $J^2(\pi)$ admits
a Hamiltonian structure induced by a vector field $X_H^1$. A vector
field $X_H^1$ on a symplectic manifold $(M,\omega)$ is called
Hamiltonian if the corresponding $\Theta = X_H^1 \rfloor \Omega $ is
closed $d\Theta =0$, and therefore exact (for a connected space).
This allows to find a Hamiltonian $H$. In this case, the symplectic
form has the form $\Omega = dx_0 \wedge dx_1$ and $\Theta = X_H^1
\rfloor \Omega =dH$. For example, for the harmonic oscillator,
\[
\textbf{H} = \int {Hdt} = \frac{{1}}{{2}}\int {\left( {x_{0}^{2} +
x_{1}^{2}} \right)dt} .
\]
Poisson structure (1)
\[
\{F,G\}=\frac{\partial F}{\partial x_1} \cdot \frac{\partial
G}{\partial x_0}-\frac{\partial F}{\partial x_0} \cdot
\frac{\partial G}{\partial x_1}
\]
is induced by the Hamiltonian vector field
\[
X_H^1=\frac{\partial H}{\partial x_0} \cdot \frac{\partial
}{\partial x_1}-\frac{\partial H}{\partial x_1} \cdot
\frac{\partial}{\partial x_0},
\]
and Cartan distribution (2) defines the dynamical Hamilton equations
for $J^2(\pi) \subset J^1( J^1(\pi))$:
\begin{equation}
\label{eq3} \overset{\cdot}{\textbf{x}}=\{H,\textbf{x}\}
\end{equation}
in such a way that the volume of Cartan differential forms (2)
\[
I = \theta _0 \wedge \theta _1 = \Omega + X_H^1\rfloor \Omega \wedge
dt
\]
is an absolute integral invariant and the presymplectic form $i
\;(di=I)$
\[
i=\frac{1}{2}(x_0 \wedge dx_1-x_1 \wedge dx_0)+H \wedge dt
\]
gives the Poincar\'e invariant of the dynamical system in question.

According to Liouville's theorem, any Hamiltonian field preserves
the volume form, i.e., the Lie derivative of the 2-form $\Omega $
with respect to the vector field $X_H^1$ is zero, $L_X \Omega = 0$.
In other words, the one-parameter group of symplectic
transformations $\{ g_t\}$ (the phase flow) generated by $X_{H}^{1}
$ leaves invariant the form $\Omega $, i.e. $g_t^\ast \Omega =
\Omega$.

\bigskip
\textbf{3.} Extending the above discussion to the case of
$J^3(\pi)$, consider the volume 3-form of the phase space
\[
\Omega=dx_0 \wedge dx_1 \wedge dx_2.
\]

\bigskip
\textbf{Theorem 1.} The volume 3-form  $ \Omega \in \Lambda ^3$
admits existence of two Hamiltonian polyvector fields $X_H^1\in
\Lambda ^1$ and $X_H^2\in \Lambda ^2$.

\bigskip
\textbf{Proof}. By definition,
\[
L_{X} \Omega = X \rfloor d\Omega + d\left( {X\rfloor \Omega} \right)
= 0.
\]
Since $ \Omega \in \Lambda ^ 3 $, we have $d\Omega = 0 $ and
\[
d\left( {X\rfloor \Omega}  \right) = 0
\]
From the Poincar\'e lemma it follows that the form $X\rfloor \Omega
$ is exact and
\[
X\rfloor \Omega = \Theta = d\textbf{H}.
\]
\noindent 1) If $X_h^1 \in \Lambda ^1$, then $ \Theta \in \Lambda ^
2 $, $\textbf{h}=h_i dx_i \in \Lambda ^ 1 $. The Hamiltonian vector
field appears as
\begin{eqnarray}
X_h^1&=&[\text{rot}\, h]_i \cdot \frac{\partial}{\partial
x_i}\label{eq4}\\
&=& \left(\frac{\partial h_2}{\partial x_1}- \frac{\partial
h_1}{\partial x_2} \right)\frac{\partial }{\partial x_0} + \left(
\frac{\partial h_0}{\partial x_2} - \frac{\partial h_2}{\partial
x_0} \right)\frac{\partial}{\partial x_1} + \left( \frac{\partial
h_1}{\partial x_0}-\frac{\partial h_0}{\partial x_1}
\right)\frac{\partial }{\partial x_2}.\notag
\end{eqnarray}
Poisson bracket (1) can be written in the form
\[
\{ \textbf{h},G\}= \left(\frac{\partial h_2}{\partial x_1}-
\frac{\partial h_1}{\partial x_2} \right)\frac{\partial G}{\partial
x_0} + \left( \frac{\partial h_0}{\partial x_2} - \frac{\partial
h_2}{\partial x_0} \right)\frac{\partial G}{\partial x_1} + \left(
\frac{\partial h_1}{\partial x_0}-\frac{\partial h_0}{\partial x_1}
\right)\frac{\partial G}{\partial x_2}.
\]
The dynamical equations have the standard form (3):
\[
\overset{\cdot}{\textbf{x}}=\{\textbf{h,x}\}
\]
2) If $X_H^2 \in \Lambda ^2$, then $ \Theta \in \Lambda ^1 $, $H \in
\Lambda ^0 $ and we obtain the Hamiltonian bivector field
\begin{equation}
\label{eq5} X_{H}^{2} = \frac{{1}}{{2}}\left( {\frac{{\partial
H}}{{\partial x_{0}} } \cdot \frac{{\partial} }{{\partial x_{1}} }
\wedge \frac{{\partial }}{{\partial x_{2}} } + \frac{{\partial
H}}{{\partial x_{1}} } \cdot \frac{{\partial} }{{\partial x_{2}} }
\wedge \frac{{\partial} }{{\partial x_{0}} } + \frac{{\partial
H}}{{\partial x_{2}} } \cdot \frac{{\partial }}{{\partial x_{0}} }
\wedge \frac{{\partial} }{{\partial x_{1}} }} \right).
\end{equation}
The Poisson bracket has more complicated expression
\begin{eqnarray}
&&X_{H}^{2} \rfloor \left( {dF \wedge dG} \right) = \left\{ H,F,G
\right\}=\frac{1}{2}\left[  {\frac{{\partial H}}{{\partial x_{0}} }
\cdot \left( {\frac{{\partial F}}{{\partial x_{1}} }\frac{{\partial
G}}{{\partial x_{2} }} - \frac{{\partial F}}{{\partial x_{2}}
}\frac{{\partial G}}{{\partial x_{1}} }} \right)}\right. \notag \\
\label{eq6}\\
&+&\left.  {\frac{{\partial H}}{{\partial x_{1}} } \cdot \left(
{\frac{{\partial F}}{{\partial x_{2}} }\frac{{\partial G}}{{\partial
x_{0} }} - \frac{{\partial F}}{{\partial x_{0}} }\frac{{\partial
G}}{{\partial x_{2}} }} \right) + \frac{{\partial H}}{{\partial
x_{2}} } \cdot \left( {\frac{{\partial F}}{{\partial x_{0}}
}\frac{{\partial G}}{{\partial x_{1} }} - \frac{{\partial
F}}{{\partial x_{1}} }\frac{{\partial G}}{{\partial x_{0}} }}
\right)}\right] \notag
\end{eqnarray}
and requires introduction of two Hamiltonians for the dynamical
equation
\begin{equation}
\label{eq7} \overset{\cdot}{\textbf{x}} = \{H ,F, \textbf{x}\}.
\end{equation}

The use of a vector Hamiltonian allows ones to obtain a
generalization of the Poincar\'e invariants. Since the volume of the
Cartan distribution
\[
\theta  = d\textbf{x}-\{ H, \textbf{x} \} \wedge dt
\]
is the absolute integral invariant $I = \theta _0 \wedge \theta _1
\wedge \theta _2$, or
\[
I = \Omega - X^{1}\rfloor \Omega \wedge dt,
\]
the pretrisymplectic form $i \quad \left( {I = di} \right)$
\[
i = \omega - (h_i dx_i) \wedge dt
\]
gives the Poincar\'e invariant of the dynamical system in question.
Here
\[
\omega =x_0 dx_1 \wedge dx_2 + x_1 dx_0 \wedge dx_2 + x_2 dx_0
\wedge dx_1
\]
is the Kirillov symplectic form $\left( {\Omega = d\omega} \right)$.

\bigskip
\textbf{Example 1.} Consider the equations of motion of a "solid
body"
\[\left\{
\begin{array}{c}
  \dot x = y-z \\
  \dot y = z-x  \\
  \dot z = x-y
\end{array}\right.
\]
The Lax pair
\[
\dot L=[ML], \qquad L=\left(
\begin{array}{lll}
x&z&y \\
z&y&x\\
y&x&z
\end{array} \right), \qquad
M=\frac{1}{2}\left(
\begin{array}{rrr}
0&-1&1 \\
1&0&-1\\
-1&1&0
\end{array} \right)
\]
for this flow gives two invariants
\[
I_1=\text{tr}\,L=x+y+z, \qquad
I_2=\frac{1}{2}\;\text{tr}\,L^2=\frac{3}{2}(x^2+y^2+z^2),
\]
and Poisson-Nambu [3] structure (7)
\[
\overset{\cdot }{x}_{i}=\{H,F,x_{i}\}=X_{H}\rfloor dF\wedge
dx_{i}=-X_{F}\rfloor dH\wedge dx_{i},
\]
is determined by vector field (6)
\begin{eqnarray*}
X_{H} &=&z\frac{\partial }{\partial x}\wedge \frac{\partial
}{\partial y}+x \frac{\partial }{\partial y}\wedge \frac{\partial
}{\partial z}+y\frac{
\partial }{\partial z}\wedge \frac{\partial }{\partial x}, \\
X_{F} &=&\frac{\partial }{\partial x}\wedge \frac{\partial
}{\partial y}+ \frac{\partial }{\partial y}\wedge \frac{\partial
}{\partial z}+\frac{
\partial }{\partial z}\wedge \frac{\partial }{\partial x}
\end{eqnarray*}
with two Hamiltonians $H=\frac{1}{3}I_2, \; F=I_1$.

To determine the Hamiltonian vector field, we rewrite the equation
of motion of a solid body in the form
\[
\overset{\cdot}{\textbf{x}}=J\textbf{x}, \qquad J=-2M.
\]
This vector flow is called Hamiltonian if
\[
\text{div} J\textbf{x}=0.
\]
This may mean that
\[
J\textbf{x}=\text{rot} \,h.
\]
The latter expression is the requirement that some differential form
is closed $d \omega=0\;$:
\[
\omega =(y-x)dy \wedge dz+(-x+z)dz \wedge dx+(x-y)dx \wedge dy.
\]
Using the homotopy operator,we obtain
\[
\omega=d \textbf{h}, \qquad \textbf{h} =(h \cdot d x),
\]
where
\[
\quad h=\frac{1}{3}\left(
\begin{array}{l}
y^2+z^2-x(y+z) \\
z^2+x^2-y(z+x)\\
x^2+y^2-z(x+y)
\end{array} \right),\quad
dx=\left( \begin{array}{l} dx \\ dy\\ dz \end{array} \right).
\]
Now the Hamiltonian vector field takes the form
\begin{eqnarray*}
X_h^1&=&\left( \text{rot}\, \textbf{h} \cdot
\frac{\partial}{\partial
\textbf{x}}\right)\\
\\
&=&\left( \frac{\partial h_3}{\partial y}-\frac{\partial
h_2}{\partial z}
 \right)\frac{\partial}{\partial
x}+ \left( \frac{\partial h_1}{\partial z}-\frac{\partial
h_3}{\partial x}
 \right)\frac{\partial}{\partial
y} + \left( \frac{\partial h_2}{\partial x}-\frac{\partial
h_1}{\partial y}
 \right)\frac{\partial}{\partial
z}\\
\\
&=&\frac{1}{3}( y^2+z^2-x(y+z))\frac{\partial}{\partial x}\\
\\
&+& \frac{1}{3}( z^2+x^2-y(z+x))\frac{\partial}{\partial y} +
\frac{1}{3}(x^2+y^2-z(x+y))\frac{\partial}{\partial z},
\end{eqnarray*}
and the Poisson bracket $\{\textbf{h},G\}=X_h\rfloor dG$ gives
\begin{equation}
\label{eq8}\overset{\cdot }{x}_{i}=\{\textbf{h},x_{i}\}=X_{h}\rfloor
dx_{i},
\end{equation}
where $\textbf{h}$ is the vector Hamiltonian and
\[
d\textbf{h}=dH\wedge dF.
\]

\bigskip
\textbf{Example 2.}  The Ishii equation [2]
$\overset{\cdots}{x}=\overset{\cdot}{x}x \in J^3(1,1)$ can be
rewritten in the form of a flow in $J^1(1,3)$:
\[
\overset{\cdot}{\textbf{x}}=J\textbf{x}, \qquad  J=\left(
\begin{array}{lll}
0&1&0 \\
0&0&1\\
y&0&0
\end{array} \right).
\]
The Lax representation for this flow
\[
\dot L=[ML], \qquad L=\left(
\begin{array}{lll}
-x&0&1 \\
-y&0&0\\
-z&\frac{y}{2}&0
\end{array} \right), \qquad
M=\frac{1}{2}\left(
\begin{array}{rrr}
2&1&0 \\
0&2&2\\
y&x&2
\end{array} \right)
\]
gives two invariants
\[
I_2=\frac{1}{2}\;trL^2=\frac{x^2}{2}-z, \qquad
I_1=\frac{1}{3}\;trL^3=xz-\frac{y^2}{2}-\frac{x^3}{3}
\]
for Eq.(7). The vector Hamiltonian appears as
\[
\textbf{h}=\frac{1}{4}\left(
\begin{array}{l}
z^2-xy^2 \\
x^2y-2yz\\
y^2-2xz
\end{array} \right)
\]
These Hamiltonians are connected by expressions
\[
d\textbf{h}=dI_1\wedge dI_2.
\]
This means that our system admit Poisson structure with vectorial
Hamiltonian
\[
\overset{\cdot }{x}_{i}=\{\textbf{h},x_{i}\}=X_{h}\rfloor dx_{i},
\]
where
\[
X_{h}=y\frac{\partial }{\partial x}+z\frac{\partial }{\partial y}%
+xy\frac{\partial }{\partial z},
\]
and Poisson-Nambu structure
\[
\overset{\cdot }{x}_{i}=\{I_1,I_2,x_{i}\}=X_{I_1}\rfloor dI_2\wedge
dx_{i}=-X_{I_2}\rfloor dI_1\wedge dx_{i},
\]
where
\begin{eqnarray*}
X_{I_1} &=&x\frac{\partial }{\partial x}\wedge \frac{\partial
}{\partial y}+(z-x^2) \frac{\partial }{\partial y}\wedge
\frac{\partial }{\partial z}-y\frac{
\partial }{\partial z}\wedge \frac{\partial }{\partial x}, \\
X_{I_2} &=&x\frac{\partial }{\partial y}\wedge \frac{\partial
}{\partial z}- \frac{\partial }{\partial x}\wedge \frac{\partial
}{\partial y}.
\end{eqnarray*}

\bigskip
\textbf{4.} In the generalization of the above discussion to the
case $J^{n}\left( {\pi} \right)$ we will rely upon Liouville's
theorem. Namely,we will seek for polyvector fields preserving the
multisymplectic volume n-form
\[
\Omega = dx_{0} \wedge dx_{1} \wedge ... \wedge dx_{n - 1}
\]
of the phase space.

\bigskip
\textbf{Theorem 2.} The multisymplectic volume $n$-form of the phase
space $\Omega \in \Lambda ^{n}$ admits existence of $n - 1$
Hamiltonian polyvector fields $X^{k} \in \Lambda ^{k}$.

\bigskip
\textbf{Proof} of this theorem is similar to the proof of Theorem 1.
It should only be noted that if $\textbf{X}_{H}^{k} \in \Lambda
^{k}$, then $\Theta \in \Lambda ^{n - k}$, $H \in \Lambda ^{n - k -
1}$. In order for $H \in \Lambda ^{m} \quad \left( {m \geq 0}
\right)$, it is necessary that $\left( {k \leq n - 1} \right)$.

A peculiarity of this problem is that $\Omega \in \Lambda ^{n}$ ,
where $\Lambda \left( {R^{n}} \right) \supset \left( {\Lambda
^{0},\Lambda ^{1},...,\Lambda ^{n}} \right) $ is a graded exterior
algebra. If the phase space is connected, a closed form $\Theta $ is
exact. But the quotient algebra $\Lambda \left( {R^{n}}
\right)/\Lambda ^{n}$ admits existence of $\left( {n - 1} \right)$
exact forms $\Theta \in \Lambda ^{m = n - k} \quad \left( {1 \leq m
\leq n - 1} \right)$, depending on whether $X^{k} \in \Lambda ^{k}
\quad \left( {n - 1 \geq k \geq 1} \right)$ or not.

It is not diffcult to present several Hamiltonian polyvector fields
\[
X_{H}^{n - 1} = \frac{{1}}{{\left( {n - 1} \right)!}}\sum\limits_{k
= 0}^{n - 1} {\frac{{\partial H}}{{\partial x_{k}} } \cdot
\frac{{\partial }}{{\partial x_{0}} } \wedge \frac{{\partial}
}{{\partial x_{1}} } \wedge \left[ {\frac{{\partial} }{{\partial
x_{k}} }} \right] \wedge ... \wedge \frac{{\partial} }{{\partial
x_{n - 1}} }} ,
\]
\[
X_{H}^{n - 2} = \frac{{1}}{{\left( {n - 2} \right)!}}\sum\limits_{i
< k}^{n - 1} {\left( {\frac{{\partial H}}{{\partial x_{i}} } -
\frac{{\partial H}}{{\partial x_{k}} }} \right) \cdot
\frac{{\partial} }{{\partial x_{0}} } \wedge \left[
{\frac{{\partial} }{{\partial x_{i}} } \wedge \frac{{\partial
}}{{\partial x_{k}} }} \right] \wedge ... \wedge \frac{{\partial
}}{{\partial x_{n - 1}} }} ,
\]
\[
X_{H}^{n - 3} = \frac{{1}}{{\left( {n - 3} \right)!}}\sum\limits_{i
< k < l}^{n - 1} {\left( {\frac{{\partial H}}{{\partial x_{i}} } -
\frac{{\partial H}}{{\partial x_{k}} } + \frac{{\partial
H}}{{\partial x_{l}} }} \right) \cdot
 \left[ {\frac{{\partial }}{{\partial
x_{i}} } \wedge \frac{{\partial} }{{\partial x_{k}} } \wedge
\frac{{\partial} }{{\partial x_{l}} }} \right] \wedge ... \wedge
\frac{{\partial} }{{\partial x_{n - 1}} }} ,
\]

\begin{center}{\dots \qquad  \dots  \qquad \dots}\end{center}

\noindent for which
\[
\Theta ^{1} = X_{H}^{n - 1} \rfloor \Omega = dH,
\]
\[
\Theta ^{2} = X_{H}^{n - 2} \rfloor \Omega = \sum\limits_{i}^{n - 1}
{dH \wedge dx_{i}} ,
\]
\[
\Theta ^{3} = X_{H}^{n - 3} \rfloor \Omega = \sum\limits_{i < j}^{n
- 1} {dH \wedge dx_{i} \wedge dx_{j}}  ,
\]
\[
\Theta ^{4} = X_{H}^{n - 4} \rfloor \Omega = \sum\limits_{i < j <
k}^{n - 1} {dH \wedge dx_{i} \wedge dx_{j} \wedge dx_{k}}  ,
\]

\begin{center}{\dots \qquad  \dots  \qquad \dots}\end{center}
\[
\Theta ^{n - 1} = X_{H}^{1} \rfloor \Omega = \sum\limits_{i < j <
... < k}^{n - 1} {dH \wedge dx_{i} \wedge dx_{j} \wedge ... \wedge
dx_{k}}  .
\]

For any $n$-form $\Omega $, the corresponding Hamiltonian vector
field $X_{H}^{1} \in \Lambda ^{1}$ generates the Poisson structure
$X_H^1 \rfloor dx_i=\{H,x_i\}$ and the contact vector field
\[
X_H=\frac{\partial }{\partial t}+X_H^1
\]
of distribution(2).

\bigskip
\textbf{5.} The Nambu theory [3] supposes introduction $\left( {n-1}
\right)$ Hamiltonians for description of the $n$-dimensional phase
space ïðîñòðàíñòâà  $\Omega $. The use of Hamiltonian polyvector
fields leads to the following generalization of this theory.

\bigskip
\textbf{Theorem 3.} The Poisson structure of an $n$-dimensional
multisymplectic space $\Omega $ is induced by $(n-1)$ Hamiltonian
polyvector fields  $X_{H}^{k} \in \Lambda ^{k}$ each of which
requires introduction of $k$ Hamiltonians:
\[
X_{H}^{k} \rfloor \left( {dF}_{1}{\wedge dF}_{2}{\wedge }...{\wedge
dF}_{k}\right) =\left\{ {F}_{1},{F}_{2},...,{F}_{k},H\right\}.
\]

\bigskip
\textbf{Example 3.} This example can be of interest for the
description of multilevel systems of the quantum information theory.
In the capacity of carriers of quantum information,at present,
observable generators of the group SO(3) (isomorphic to SU(2)) are
considered. To obtain irreducible representations of the group
SU(2), the standard method of constructing proper states of the
angular momentum operator is applied. In addition, the state
indication operator
\[
\Sigma_3\left| {n,m} \right\rangle=m\left| {n,m} \right\rangle
\]
and the creation-destruction operators (level increasing or
decreasing operators) with boundary properties
\[
\Sigma_+\left| {n,n} \right\rangle=\Sigma_-\left| {n,-n}
\right\rangle=0
\]
are considered.

Another way ([4].P.8) is to use annular $n$-level systems with
periodicity condition
\[
\left| {n,i+n} \right\rangle=\left| {n,i} \right\rangle.
\]
This allows ones to use in the description the two operators
\[
\Sigma_1\left| {n,i} \right\rangle=\left| {n,i+1} \right\rangle (mod
\;n), \quad \Sigma_3\left| {n,i} \right\rangle=i\left| {n,i}
\right\rangle,
\]
where
\[
\Sigma_1 =\left(
\begin{array}{ccccc} 0&&&&1 \\ 1&0&&& \\ &1&\ddots &&
\\ &&\ddots&0& \\ &&&1&0 \end{array}
\right), \quad \Sigma_3 =\left(
\begin{array}{ccccc} \sigma^0&&&& \\ &\sigma^1&&& \\ &&\sigma^2&& \\ &&&\ddots &
\\ &&&&\sigma^n \end{array}
\right),
\]
Here $\sigma=e^{i\frac{2\pi}{n}}$ has the following properties:
\[
\sigma^n=1, \quad \sigma^{n+k}=\sigma^k(mod \;n),\quad
\sigma=\sigma^{n+1}, \quad \sum \limits_{k=0}^n\sigma^k=0.
\]
In what follows the generators $\Sigma_1$ and $\Sigma_3$ will be
called the generalized Pauli matrices. We will denote the identity
matrix by $\textbf{I}=\Sigma_0$ and the increasing operator by
\[
\Sigma_1^+=\Sigma_1^{n-1} =\left(
\begin{array}{ccccc} 0&1&&& \\ &0&1&& \\ &&\ddots &&
\\ &&&0&1 \\ 1&&&&0 \end{array}
\right).
\]
In particular, for $n=2$, the expansion of the operator
$e^{t\Sigma_1}$ into a series gives
\[
e^{t\Sigma_1}=c_0(t)\Sigma_0+c_1(t)\Sigma_1,
\]
where
\begin{equation}
\label{9} c_0(t)=\sum \limits_{k=0}^\infty \frac{t^{2k}}{(2k)!},
\quad c_1(t)=\sum \limits_{k=0}^\infty \frac{t^{2k+1}}{(2k+1)!}
\end{equation}
In this case, we let
\[
c_0(t)=\cosh(t), \quad c_1(t)=\sinh (t)
\]
and, taking into account that
\[
\frac{d}{dt}\cosh (t)=\sinh (t), \quad \frac{d}{dt}\sinh (t)=\cosh
(t)
\]
obtain the system of differential equations
\begin{equation}
\label{eq11}\dot {\textbf{c}}=\Sigma_1^+ \textbf{c},
\end{equation}
where $\textbf{c}(t)$ is the vector with coordinates
$\textbf{c}=(c_0,c_1)$. The Lax representation appears as
\[
L=\left(
\begin{array}{ll}
c_1y&c_0 \\ c_0\sigma ^1&c_1\sigma ^1
\end{array} \right), \qquad
M=\frac{\sigma ^1}{2}\;\Sigma_1
\]
It gives the scalar invariant
\[
I=\frac{1}{2}\;trL^2=c_1^2-c_0^2
\]
for Poisson bracket (3). Note that the relation $I=1$ is the basic
trigonometric identity for functions (9).

For $n=3$,the expansion of the operator $e^{t\Sigma_1}$ into a
series gives
\[
e^{t\Sigma_1}=c_0(t)\Sigma_0+c_1(t)\Sigma_1+c_2(t)\Sigma_1^2,
\]
where
\begin{equation}
\label{eq11}c_0(t)=\sum \limits_{k=0}^\infty \frac{t^{3k}}{(3k)!},
\quad c_1(t)=\sum \limits_{k=0}^\infty \frac{t^{3k+1}}{(3k+1)!},
\quad c_2(t)=\sum \limits_{k=0}^\infty \frac{t^{3k+2}}{(3k+2)!}.
\end{equation}
The obtained coeffcients again satisfy system of differential
equations (10), which has the following form in the Lax
representation:
\[
L=\left(
\begin{array}{lll}
c_0\sigma ^2&c_1\sigma ^1&c_2\sigma ^0 \\
c_2\sigma ^1&c_0\sigma ^0&c_1\sigma ^2\\
c_1\sigma ^0&c_2\sigma ^2&c_0\sigma ^1
\end{array} \right), \qquad
M=\frac{1}{\sigma ^1-\sigma ^2}\;\Sigma_1.
\]
This system gives the only one scalar invariant
\[
I=\frac{1}{3}\;trL^3=c_0^3+c_1^3+c_2^3-3c_0c_1c_2.
\]
Note that the expression $I=1$ is the basic trigonometric identity
for functions (11). The vector Hamiltonian for bracket (8) appears
as
\[
\textbf{h}=\frac{1}{4}\left(
\begin{array}{l}
c_2^2-2c_0c_1 \\
c_0^2-2c_1c_2 \\
c_1^2-2c_0c_2
\end{array} \right).
\]

In the general case, we have
\[
e^{t\Sigma_1}=c_0(t)\Sigma_0+c_1(t)\Sigma_1+...+c_{n-2}(t)\Sigma_1^{n-2}+c_{n-1}(t)\Sigma_1^{n-1},
\]
where the functions
\[
c_j(t)=\sum \limits_{k=0}^\infty \frac{t^{nk+j}}{(nk+j)!}
\]
satisfy system (10) with Lax pair
\[
L=\left(
\begin{array}{ccccccc}
c_{0}\sigma ^{n-1} & c_{1}\sigma ^{n-2} & c_{2}\sigma ^{n-3} & \dots
&
c_{n-3}\sigma ^{2} & c_{n-2}\sigma ^{1} & c_{n-1}\sigma ^{0} \\
c_{1}\sigma ^{n-2} & c_{0}\sigma ^{n-3} & c_{1}\sigma ^{n-4} & \dots
&
c_{n-4}\sigma ^{1} & c_{n-3}\sigma ^{0} & c_{n-2}\sigma ^{n-1} \\
c_{2}\sigma ^{n-3} & c_{3}\sigma ^{n-4} & c_{0}\sigma ^{n-5} & \dots
&
c_{n-5}\sigma ^{0} & c_{n-4}\sigma ^{n-1} & c_{n-3}\sigma ^{n-2} \\
\dots  & \dots  & \dots  & \dots  & \dots  & \dots  & \dots  \\
c_{3}\sigma ^{2} & c_{4}\sigma ^{1} & c_{n-1}\sigma ^{0} & \dots  &
c_{0}\sigma ^{5} & c_{1}\sigma ^{4} & c_{2}\sigma ^{3} \\
c_{2}\sigma ^{1} & c_{3}\sigma ^{0} & c_{4}\sigma ^{n-1} & \dots  &
c_{n-1}\sigma ^{4} & c_{0}\sigma ^{3} & c_{1}\sigma ^{2} \\
c_{1}\sigma ^{0} & c_{2}\sigma ^{n-1} & c_{3}\sigma ^{n-2} & \dots &
c_{n-2}\sigma ^{3} & c_{n-1}\sigma ^{2} & c_{0}\sigma ^{1}%
\end{array}%
\right),
\]
\[
M=\frac{1}{\sigma ^1-\sigma ^{n-1}}\;\Sigma_1.
\]
The vector invariant has a rather complicated form. For example, for
$n=4$, we have the scalar invariant
\[
I=2(2c_0c_2-c_1^2-c_4^2),
\]
and the differential form of vector Hamiltonian
\begin{eqnarray*}
\textbf{h}=(h \cdot dc)&=&\frac{1}{6}\left[
(c_3^2-2c_0c_2)dc_0 \wedge dc_1+ (c_1^2-2c_2c_0)dc_2 \wedge dc_3 \right.\\
&+&(c_2^2-2c_3c_1)dc_0 \wedge dc_3
+(c_0^2-2c_3c_1)dc_2 \wedge dc_1 \\
&+&\left.2(c_3c_0-c_1c_2)dc_1 \wedge dc_3+ 2(c_0c_1 -c_2c_3)dc_0
\wedge dc_2\right].
\end{eqnarray*}
In all the cases, the invariants are in involution, and the Poisson
bracket $\{ \textbf{h},I\}=0$ gives no new integrals.

\begin{center}
References
\end{center}

\begin{enumerate}
  \item P.A. Griffiths, {\it Exterior differential systems and the
calculus of variations}, Birkhauser, Boston (1983).
  \item M.Ishii, Painleve Property and Algebraic Integrability of Single
Variable Ordinary Differential Equations with Dominants,
// Progr.Theor.Phys.84,386–391(1990).
  \item Y. Nambu, Generalized Hamiltonian dynamics, {\it Phys.Rev.D}, {\bf
7} (1973), 5405-5412.
  \item K. Fujii, A New Algebraic Structure of Finite Quantum Systems and
the Modified Bessel Functions, arXiv:0704.1844.
\end{enumerate}

\bigskip
Translated by V.V.Shurygin

\end{document}